
\input amssym.def
\input amssym

\magnification = \magstep1
\hsize=14.4truecm
\vsize= 22truecm
\hoffset=0.4 truecm
\voffset=0.3 true in
\baselineskip=12pt

\font\big=cmbx10 scaled \magstep2
\font\caps=cmcsc10
\font\fascio=eusm10
\font\titol=cmcsc8
\font\eightpt=cmr8

\def\smallskip{\vskip 0.2truecm}
\def\medskip{\vskip 0.5truecm}
\def\bigskip{\vskip 0.7truecm}
\def\tit{\medskip \noindent} 
\def\proof{\par\smallskip\noindent{\it Proof.}\ }
\def\endproof{\hfill$\square$\smallskip}

\def\sp#1{\Bbb P^{#1}}
\def\ga{p_{\scriptscriptstyle a}}

\headline={\ifnum\pageno=1\hfil\else
{\ifodd\pageno
\rightheadline\else\leftheadline\fi}\fi}
\def\rightheadline{\titol\hfil 
Linear systems representing scrolls on a rational surface
\hfil{\tenrm\folio}}                    
\def\leftheadline{{\tenrm\folio}\titol\hfil 
Mezzetti - Portelli \hfil}

\nopagenumbers

{\big 
\hbox{Linear systems representing threefolds}}

\medskip

{\big 
\hbox{which are scrolls on a rational surface
}}

\bigskip

\noindent{\it Emilia Mezzetti \ } and {\ \it 
Dario Portelli}

\medskip

{\eightpt {\noindent Dipartimento di Scienze 
Matematiche, \ Universit\`a di Trieste,

\noindent Piazzale Europa 1,
 
\noindent 34127 Trieste, ITALY

\noindent e-mail: mezzette@univ.trieste.it, 
porteda@univ.trieste.it}}

\bigskip\bigskip
\qquad\qquad\qquad\qquad\qquad\qquad {\it Dedicated to the memory of Ferran Serrano}
\bigskip\bigskip\bigskip


\tit
{\bf Abstract}
\tit


Let $X$ be a scroll over a rational surface. We
construct a linear system of surfaces in $\sp 3$
yielding a birational map $\sp 3\dasharrow X.$
We apply this construction to the scrolls of
Bordiga and Palatini.

\bigskip


\tit
{\bf Introduction}
\tit


Let $X\subset\sp r$ be a normal threefold which is a 
scroll over a rational surface. Then $X$ is a rational 
variety. 
In this paper we will expose a method to determine
a $r$-dimensional linear system $\vert\Sigma\vert$ 
of surfaces in $\sp 3$ such that the rational map 
$f:\sp 3\dashrightarrow\sp r$ associated to 
$\vert\Sigma\vert$ maps $\sp 3$ birationally onto $X.$

Our interest for this problem arose when studying 
smooth threefolds of $\sp 5.$
In [4] we considered smooth rational threefolds 
of $\sp 5$ with rational hyperplane section: it is known 
that there are exactly five families of such varieties, 
of degrees $d$, with $3\leq d\leq 7.$ 
Precisely, for $d=3$ there is $\sp 2\times\sp 1$, 
for $d=4$ the Del Pezzo threefold, complete intersection 
of two hyperquadrics, for $d=5$ the Castelnuovo threefold, 
for $d=6$ the Bordiga scroll, for $d=7$ the Palatini scroll 
(for more reference, see [7]). 
We proved that for $3\leq d\leq 6$ for any such threefold
$X$ there exists a line $L\subset X$ such that the projection
centered at $L$ from $X$ to $\sp 3$ is birational. 
In this way, we obtained a \lq\lq\thinspace uniform'' 
description of the wanted linear systems $\vert\Sigma\vert.$
Conversely, if $d=7,$ i.e. if $X$ is a Palatini scroll, 
we proved that a line with the above property never exists.

It remained the problem of finding a linear system 
$\vert\Sigma\vert$ in this last case. The method we have 
used for solving the problem is based on the fact that 
Palatini threefold is a scroll over a rational surface, 
precisely over a smooth cubic surface of $\sp 3.$ 
In fact a similar construction can be made for any 
scroll $X$ over a rational surface $V.$
The basic idea is classical: it is quoted e.g. in the 
paper [3] by Jongmans. Roughly speaking, it 
consists in establishing first an explicit birational 
correspondence between the lines of the scroll $X$ 
and the lines of $\sp 3$ passing through a fixed point 
$P$, and then a map between any two corresponding lines
is given by projection from a fixed center. 

More precisely: assume that $X$ is contained in $\sp r$ 
and fix a general linear subvariety $\Lambda$ of $\sp r$ of 
codimension $2$. 
We compose the scroll map $X\to V$ with a fixed birational 
map $V\dashrightarrow \sp 2$ to get a rational map 
$h:X \dashrightarrow \sp 2.$ 
Then we have an explicit birational map $\alpha :X
\dashrightarrow\sp 2\times\sp 1$ defined by $\alpha =
(h,\pi _{\Lambda})$, where $\pi_{\Lambda}:X\dashrightarrow
\sp 1$ is the projection centered at $\Lambda.$ 
Finally, $\alpha$ has to be composed with the projection 
$\pi _L: \sp 2\times\sp 1 \dashrightarrow\sp 3$
from a line $L$, contained in one of the planes of 
$\sp 2\times\sp 1$. Note that the projection $\pi _L$ 
which is clearly birational, contracts such 
plane to a point $P$.

Here we perform a detailed study of such birational map 
$X\dashrightarrow\sp 3$ and of its birational inverse from 
a modern point of view. The fundamental step is the complete
determination of the exceptional divisors and the fundamental
loci for both maps. From this it follows readily a rather
satisfactory description of the linear system $\vert\Sigma\vert$
of surfaces of $\sp 3$ defining the map $\sp 3\dashrightarrow X,$
in particular of its base locus. Experience shows the importance
of the knowledge of the {\sl arithmetic} genus of the $1$-dimensional
components of $Bs(\vert\Sigma\vert )$; we show how this can be done
by working-out the case of the Bordiga scroll.

It is interesting to note that the birational maps from
$\sp 3$ to the Bordiga scroll obtained in
this paper and in [4] are different. Therefore, by composition,
we obtain elements of the Cremona group
of $\sp 3$, given in general by surfaces
of degree $5$ and $7$. 
For special choices made in the construction
of the maps, we get classical cubo--cubic
transformations.

The content of the various sections is as follows.
Some basic facts concerning birational maps are recalled
for the reader's convenience in \S\thinspace 1.
In \S\thinspace 2 we study the birational map 
$g:X\dashrightarrow\sp 3$ described above; in particular 
we find its exceptional divisors and the exceptional 
divisors of the inverse of $g.$
\S\thinspace 3 is devoted to the study of the linear system
$\vert\Sigma\vert ,$ and contains
the main result of the paper: Theorem 3.1. 
The statement of this theorem
concerns the determination of the irreducible 
components of $Bs(\vert\Sigma\vert ),$ their
mutual positions and their infinitesimal structure.
Finally, in \S\thinspace 4 we apply Theorem 3.1 
to the scrolls of Bordiga and Palatini.
The results on $\vert\Sigma\vert$ obtained in these
two cases are collected in the following table (where we
denote by $A\cdot B$ the length of the zero dimensional
scheme $A\cap B$):
 
\bigskip

\halign{
\hfil\bf# & \quad\hfil#\hfil & \quad\hfil#\hfil & \quad#\hfil                        
\cr
Variety&{\bf Degree}&$deg\thinspace\Sigma$&{\bf Base locus}                           
\cr
&&&\cr
Bordiga scroll     & 6 & 5 &  
$\bullet$\ a curve $B_2$ of degree $14,$ \  geometric\cr
&&&\ \ \thinspace genus $3,$ arithmetic genus $23,$\ with\cr
&&&\ \ \thinspace an ordinary multiple point $P$ of\cr
&&&\ \ \thinspace multiplicity $10$\cr
&&&  $\bullet$\ a line $B_1$ not containing $P,$  with\cr
&&&\ \ \thinspace $B_1\cdot B_2=4$\cr
&&&\cr                            
Palatini scroll    & 7 & 7 &  
$\bullet$\ a curve $B_2$ of degree $11,$ \  geometric\cr
&&&\ \ \thinspace genus $4,$ arithmetic genus $10,$ \  with\cr
&&&\ \ \thinspace an ordinary multiple point $P$ of\cr
&&&\ \ \thinspace multiplicity $5$  \cr
&&&  $\bullet$\ a line $B_1$ not containing $P,$  with\cr
&&&\ \ \thinspace $B_1\cdot B_2=6$  \cr
&&&  $\bullet$\ the first infinitesimal neighborhoods\cr
&&&\ \ \thinspace of six lines \ $B_i$, \ each containing \ $P,$\cr
&&&\ \ \thinspace with $B_1\cdot B_i=0$ and $B_2\cdot B_i=2$\cr}

\bigskip

We wish to thank Silvio Greco
for some useful discussion. 

\smallskip

This work has been done in the framework of the 
activities of Europroj. Both authors have been
supported by funds of MURST and  
GNSAGA, Progetto Strategico \lq\lq\thinspace
Applicazioni della Matematica per la Tecnologia
e la Societ\` a", Sottoprogetto \lq\lq\thinspace
Calcolo simbolico".


\tit
{\bf  1.-\ Some property of birational maps.}
\tit


We will always work over an algebraically closed field of
characteristic $0$.
In this section we collect some classically known 
facts about birational maps. Moreover, we give the 
main features of a particular birational map, namely
a projection $\sp 2\times\sp 1\dashrightarrow\sp 3,$ 
which we shall need in the subsequent sections. 
For precise references to the literature and for proofs  
see [4].

\smallskip

Let $f:Y\dashrightarrow X$ denote 
a birational map between {\sl normal} varieties. 
A point $y\in Y$ is called {\sl fundamental} for 
$f$ if $f$ is not regular at $y.$ 
We will denote by $F$ the set of the fundamental 
points for $f,$ and we will call it the {\sl 
fundamental locus} for $f.$ 

The importance of the fundamental locus for our 
investigation is due to the fact that, {\sl if a birational
map $f:\sp n\dashrightarrow X$ is defined by a linear 
system $\vert\Sigma\vert$ of hypersurfaces in $\sp n,$
then the irreducible components of the fundamental 
locus of $f$ are precisely those of the base locus 
$Bs(\vert\Sigma\vert )$ of $\vert\Sigma\vert .$}

\medskip 

To determine the fundamental locus the starting 
point is

\smallskip

\noindent
{\bf Van der Waerden'\thinspace s Purity Theorem.}\ 
\ {\sl  With the notations introduced above,
let $W\subset Y$ be an irreducible component of the 
fundamental locus of $f.$ Let us denote by $g$ the
birational inverse of $f$ (we will mantain this notation).
Assume that $g(X)\cap W$ 
is dense in $W$ and  that $W$ is not contained in the 
singular locus of $Y.$  Then any component $E$ of 
$g^{-1}(W)$ is of codimension $1$ in $X.$} 

\smallskip

We will call any $E\subset X$ as above an 
{\sl exceptional divisor} for $g.$ Since 
the fundamental locus $F$ for $f$ is a closed 
subset of $Y$ and $codim_Y (F)\geq 2,$ every
irreducible component of $F$ is of the form 
$g(E)$ for some exceptional divisor $E$ of $g.$ 
Therefore, once we know the exceptional divisors for 
$g$, we will be able to describe $Bs(\vert\Sigma\vert ).$

\medskip

Let us recall that a {\sl characteristic
curve} $\Gamma$ of a linear system 
$\vert\Sigma\vert$ of surfaces of $\sp 3$
is the free intersection of two general surfaces
of $\vert\Sigma\vert$.
If $E$ is an exceptional surface for $f$, 
the birational map associated to 
$\vert\Sigma\vert ,$ and $B$ is the corresponding
base curve, then $deg \ E=\# (B\cap\Gamma)$.

\medskip

In the next section we will need the following
example.

\smallskip

Let $X\subset\sp 5$ be the image of $\sp 2\times\sp 1$
embedded in $\sp 5$ by the Segre map $s$. 
For a fixed point $a\in\sp 1,$ we will denote 
$F_1:= s(\sp 2\times a),$ a plane on $X.$
Moreover, let $l$ be a line in $\sp 2$ and
set $F_2:=s(l\times\sp 1).$ 
$F_2$ is a quadric surface on $X.$ Finally, 
denote by $L$ the line $F_1\cap F_2.$  Then:

\goodbreak 

\noindent
{\bf Proposition 1.1} {\sl
\noindent\item{(i)} 
The projection map $ 
\pi_L :X \dashrightarrow  \sp 3 $ 
with center $L$ is birational;

\noindent
\item{(ii)} the exceptional divisors of $\pi_L$ are
$F_1$ and $F_2,$ in particular $\pi_L(F_1)$ is a 
single point $P$ of $\sp 3,$ and $\pi_L(F_2)$ 
is a line $B\subset\sp 3$ such that $P\notin B;$

\noindent\item{(iii)} 
the map $\pi_L^{-1}$ is defined by the 
linear system  of quadrics $\vert\Sigma\vert$ with 
base locus $B\cup P;$

\noindent\item{(iv)} 
the only exceptional divisor for $\pi_L^{-1}$ 
is the plane $\Phi =\langle B\cup P\rangle ,$ which is 
contracted to the line $L .$ 

}


\tit
{\bf  2.-\ Construction of a birational map from $\sp 3$ to 
$X$.}
\tit


Let $X\subset\sp r$ be a normal scroll of degree $d$
over a rational  surface $V.$ 
Let $u:X\to V$ denote the scroll map. 
We will construct, now, a birational map 
$g:X\dashrightarrow\sp 3.$

\smallskip

Let us fix a birational map $v:V\dashrightarrow\sp 2$ 
and set $h:=v\circ u:X\dashrightarrow\sp 2.$ 
Moreover, we fix inside $\sp r$ a linear subvariety 
$\Lambda\simeq\sp {r-2}$ such that $\Lambda$ does not 
contain any line of the scroll $X.$ Let
$\pi_\Lambda :\sp r\dashrightarrow\sp 1$ denote the 
projection map. Then we have a rational map

$$
\alpha :X\dashrightarrow\sp 2\times\sp 1$$
defined by
$$\alpha (x):=(h(x),\pi_\Lambda (x))
$$
which clearly results to be birational 
(see Figure 1 below).

Let $L\subset\sp 2\times\sp 1$ be a line  
as in \S\thinspace 1. 
If we compose $\alpha$ with the projection $\pi_L$ we 
get a birational map \ $g:=\pi_L\circ\alpha 
:X\dashrightarrow\sp 3.$ 

Let us remark that the lines of the scroll $X$
correspond via the map $g$ to the lines in $\sp 3$
through the fixed point $P=\pi_L(F_1).$

\smallskip

Let $f$ be the birational inverse of $g$ and let
 $\vert\Sigma\vert$ be 
the linear system of surfaces in $\sp 3$ defining 
 $f.$ 
The main goal of this section is determination of
all the exceptional divisors for $g.$  

Since $g$ is the composition of the birational maps 
$\alpha$ and $\pi_L,$ to find out its exceptional 
divisors  it is sufficient to determine the 
exceptional divisors for $\alpha$ and $\pi_L,$ separately. 
For $\pi_L$ this was already done in Proposition 1.1. 

The exceptional divisors for $\alpha$ are of two 
different kinds, i.e. coming either from $h$ or from 
$\pi _{\Lambda}.$
First of all, if $\Delta\subset V$ is an exceptional 
divisor for the map $v:V\dashrightarrow\sp 2,$ then 
$\alpha$ contracts the surface $E:=u^{-1}(\Delta )$ 
to the line $v(\Delta )\times\sp 1\subset
\sp 2\times\sp 1$ (see Fig.\thinspace 1).

\null
\vskip 8cm

\centerline{-- Fig.\thinspace 1 --}

\smallskip

On the other hand,  $\alpha$ is not regular on 
the curve $C:=X\cap\Lambda .$ Nevertheless, it is clear that 
$\alpha$ is constant on any line  of the scroll $X$ which 
meets $C;$ therefore $\alpha$ contracts the ruled surface 
$E_{\scriptscriptstyle {\infty}}:=u^{-1}(u(C))$ to a curve 
$C'\subset\sp 2\times\sp 1.$ Clearly, there is no other 
exceptional divisor for $\alpha ,$ and we conclude:

\smallskip

\noindent
{\bf Proposition 2.1} {\sl Let $F_1$ and $F_2$ be the 
surfaces on $\sp 2\times\sp 1$ defined in Proposition 1.1, 
and set $E_i:=\alpha^{-1}(F_i)$ where $i=1,2.$ Let 
$\Delta_{\scriptscriptstyle 3},\ldots ,\Delta_s$ be the 
exceptional divisors for $v$ on $V,$ and set 
$E_i :=u^{-1}(\Delta_i),$ for $i=3,\ldots ,s.$
The exceptional divisors for $g$ are the surfaces 
$E_{\scriptscriptstyle 1},E_{\scriptscriptstyle 2},
E_{\scriptscriptstyle 3}\ldots ,E_s,
E_{\scriptscriptstyle {\infty}}.$}

\smallskip

\noindent
Note, in particular, that $E_1$ is a hyperplane 
section of $X$ by construction of $\alpha .$

\medskip 

We can determine the 
exceptional divisors for the map $f:\sp 3\dashrightarrow X$
by a similar \lq\lq\thinspace step by step" procedure, by
determining the exceptional divisors for $\pi _L^{-1}$ and for
$\alpha ^{-1}$ separately.

Recall from Proposition 1.1 that 
the unique exceptional divisor for
$\pi _L^{-1}$ is the plane generated by the point 
$P:=\pi_L(F_1)$ 
and the line $B_1:=\pi_L(F_2).$ 

As for the exceptional divisors for $\alpha ^{-1}$ they
arise either from $h$ or from $\pi _{\Lambda}.$
Let $T\subset\sp 2$ be an exceptional divisor 
of $v^{-1}:\sp 2\dashrightarrow V.$ Then 
$\alpha ^{-1}$
contracts the surface $T\times\sp 1$ to the line
$u^{-1}(v^{-1}(T))=h^{-1}(T)$.

There is one more exceptional divisor, arising from 
$\pi _\Lambda .$ Let $\Gamma :=h(C)\subset\sp 2$ and  set
$U:=\Gamma\times\sp 1\subset\sp 2\times\sp 1$
(note that the curve $C'$ 
defined above lies on $U;$ see Fig.\thinspace 2).
Finally, set $\Phi :=\pi _L(U).$

\medskip

\noindent
{\bf Proposition 2.2} {\sl  
The exceptional divisors for $f$ are the plane 
$\langle P\cup B_1\rangle$, the surfaces
$\pi_L(T\times\sp 1),$ where $T$ is as above, and the 
surface $\Phi.$ Moreover, $\Phi$ is a cone 
with vertex $P$ and degree $n,$ where $n=deg 
\thinspace\Gamma .$}

\proof 
For any $b\in\Gamma$ the line $b\times\sp 1$ is contained
in $U,$ and it is mapped by $\pi_L$ to a line through $P.$ 
Hence $\Phi =\pi_L(U)$ is a cone with vertex $P.$ 
It remains to show that $\Phi$ is a cone of degree $n.$
A line in $\sp 2$ cuts $\Gamma$ in $n$ points. Therefore, a 
line $R\subset\sp 2\times a$ cuts $U$ in $n$ points. 
Since $\pi_L(R)$ is a line in $\sp 3,$
the degree of the cone $\Phi$ is $n.$

To prove that $\Phi$ is an exceptional divisors for $f$ 
it is sufficient to show that $\alpha ^{-1}$ contracts 
the surface $U$ to the curve $C\subset X$ defined above. 
More precisely it is sufficient to show that $\alpha ^{-1}$ 
contracts any line $Q\times\sp 1$, where $Q\in\Gamma$, to a 
point of $C.$ 
In fact, if $Q'=(Q,a)\in Q\times\sp 1$ is a general 
(=\lq\lq\thinspace not in $C'$\thinspace\thinspace ") 
point of the line, then 
$\alpha ^{-1}(Q')$ is the intersection of the line
$h^{-1}(Q')$ of the scroll with a suitable hyperplane 
$H\subset\sp r$ containing $\Lambda .$ But this intersection
is clearly a point on $C$ independent of the hyperplane $H.$
\endproof

\medskip

For further reference we make the following

\smallskip 

\noindent
{\bf Remark 2.3} The planes
$\sp 2\times a$ on $\sp 2\times\sp 1$ correspond    
via the map $\pi_L$ to the planes $W_a$ in $\sp 3$
containing the line $B_1$  and 
via the map $\alpha$ to the sections $S_a$ of $X$ 
with the hyperplanes $H_a$ containing $\Lambda .$


\tit
{\bf  3.-\ Description of the linear system defining 
$f:\sp 3\dashrightarrow X$.}
\tit


The construction of the birational maps $f$ and $g$ made 
in the previous section depends on several choices.
We will make now some harmless \lq\lq transversality" 
assumptions; under them we will get a very precise 
description of $Bs(\vert\Sigma\vert ).$

First of all, we can assume $C\cap Sing(X)=\emptyset .$
In fact, $X$ is normal, hence $codim_XSing(X)\geq 2.$

Now, let us consider the rational surfaces 
$E_i$ for $i=3, \ldots ,s$ (see Prop. 2.1). Assume 
that $E_i$ has degree $\delta_i$ in $\sp r.$
In particular, being $C$ a curve section of $X,$ we have 
$C\cdot E_i=\delta_i .$ 
By moving $\Lambda ,$ hence $C,$ it is harmless to assume 
that $\bigcup_i\thinspace\thinspace u(C\cap E_i)$ are 
$\sum_i\delta_i$ \ distinct points on $V.$

Moreover, we will assume that the hyperplane $H\subset\sp r$
such that $E_1=X\cap H$ does not contain any line of the 
scroll $X$ through a point of $C\cap E_i,$ for any $i=3, 
\ldots ,s.$

\smallskip
We can state now our main result:

\smallskip

\noindent
{\bf Theorem 3.1} {\sl   

Let $X\subset \sp r$ be a normal scroll of degree $d$ 
over a rational surface. Let $f$ be the birational map
from $\sp 3$ to $X$ defined in \S 2, and let $\vert\Sigma\vert$
be the linear system of surfaces defining $f$. 

Keeping the
notations introduced in \S 2,
under the transversality assumptions made at the beginning
of this section, we have:

\noindent\item{(i)} The $1$-dimensional components of 
$Bs(\vert\Sigma\vert )$ are the line $B_{\scriptscriptstyle 1}
:=\pi_L(F_2),$ the lines $B_i:=\pi_L(v(\Delta_i)\times\sp 1)$ 
for $i=3,\ldots ,s$ 
and the curve $B_{\scriptscriptstyle 2}:=\pi_L(C')\subset\Phi .$  
The unique $0$-dimensional component of $Bs(\vert\Sigma\vert )$ 
is the point $P=\pi_L(F_1)$. 

\noindent\item{(ii)} The curve $B_2$ intersects the line $B_1$ 
at $n$ points (where $n=deg (\Gamma )$) and is $\delta_i$-secant to any line 
$B_i,$ \ $i=3, \ldots ,s.$ \ 
$B_2$ is unisecant all the other lines on $\Phi .$

\noindent\item{(iii)} The degree of $B_2$ is equal to
$deg \ E_{\scriptscriptstyle {\infty}} -d+n \ $ 
and its geometric genus is equal to the sectional 
genus of $X.$ The point $P$ is a multiple point of 
multiplicity $deg \ E_{\scriptscriptstyle {\infty}} 
- d$ for $B_2.$
The tangent lines of $B_2$ at $P$ are distinct and
the branches of $B_2$ at $P$ are nonsingular (for the
terminology concerning singular points we 
follow {\rm [5]}; see also \S\thinspace 4). 
If the sectional genus of $X$ is strictly positive, 
then $B_2$ is smooth outside $P$. 

\noindent\item{(iv)} The surfaces $\Sigma$ are monoids 
of degree $n+1$ whose point of multiplicity $n$ is $P.$ 
The scheme structure of the $0$-dimensional (embedded)
component $P$ of the 
base locus $Bs(\vert\Sigma\vert )$ is exactly that of 
$(n-1)$-th infinitesimal neighborhood of $P.$ Let $S$
be a hyperplane section of $X$ corresponding to
$\Sigma\in\vert\Sigma\vert ,$ and set $C_1:=S
\cap E_1.$ Then the tangent cone of $\Sigma$ at $P$
is $\pi_L(h(C_1)\times\sp 1).$

\noindent\item{(v)} The base
locus of $\vert\Sigma\vert $ contains the $(\delta_i-1)$-th
infinitesimal neighborhood of $B_i,$ for any $i=3, \ldots ,s.$ 

}

\proof (i) follows from Propositions 2.1 and 1.1. 
In particular, by definition, $B_2=
g(E_{\scriptscriptstyle {\infty}})=\pi_L(C').$ 
As remarked above, we have $C'\subset U,$ hence 
$B_2$ lies on $\Phi .$

On $\sp 2\times\sp 1$ we have $F_2=L\times\sp 1.$
Then, from $deg\ \Gamma =n$ it follows that 
$C'\cdot F_2=n$ and $B_2$ intersects the line $B_1$ 
at $n$ points. To see that $B_2$ is $\delta_i$-secant 
to $B_i,$ it is 
sufficient to remark that all the $\delta_i$ points of 
$C\cap E_i$ have the same image in the map $h :
X\dashrightarrow\sp 2,$ by definition of $E_i$ 
(see Fig. 2, below).

The fact that $B_2$ is unisecant to any line on the 
cone $\Phi$ which is different from $B_3,\ldots ,B_s$
is easily checked, and this completes the proof of (ii).

To compute the degree of $B_2,$ let us consider a 
plane $W\subset\sp 3$ such that $B_1\subset W.$
By Remark 2.3, the points of
$W\cap B_2$ outside $B_1$ correspond to points of 
$(\sp 2\times a)\cap C',$ hence they correspond to 
lines on $S_a\cap E_{\scriptscriptstyle {\infty}}.$
Since $C\subset H_a,$ from the definition of
$E_{\scriptscriptstyle {\infty}}$ we get 
$H_a\cap E_{\scriptscriptstyle {\infty}}
=C\cup l_1\cup\ldots\cup l_p,$ where every 
$l_j$ is a line of the ruling. By Bezout: 
$p=deg \ E_{\scriptscriptstyle {\infty}}-d.$ 
Since $B_1$ meets $B_2$ at $n$ points, we get
the announced degree for $B_2.$

\null
\vskip 8 cm

\centerline{-- Fig.\thinspace 2 --}

\smallskip

From the discussion above it follows that any plane 
$\sp 2\times a$ on $\sp 2\times\sp 1$ cuts $C'$ 
at $deg\ E_{\scriptscriptstyle {\infty}}-d$ points.
In particular, taking $\sp 2\times a=F_1$ we conclude that
$P$ is a multiple point of $B_2$ of multiplicity 
$deg\thinspace E_{\scriptscriptstyle {\infty}}-d.$

To study the singular point $P$ of $B_2$ we start
by considering the common points of $C'$ and $F_1.$ 
These points are in one-to-one correspondence with the
lines on $E_{\scriptscriptstyle {\infty}}$ contained in
the hyperplane $H\subset\sp r$ such that $E_1=X\cap H.$
We let $P_1,\ldots ,P_r\in\sp 2$ denote the image in $h$
of such lines.
Since the curve $C$ is smooth and any line of the scroll 
$E_{\scriptscriptstyle {\infty}}$ intersects $C$ in 
exactly one point, {\sl the curve $C'$ is smooth.} 
The singular locus of the ruled surface 
$U\subset \sp 2\times\sp 1$ is formed by the lines 
$b\times\sp 1,$ where $b$ is a singular point of the 
curve $\Gamma =h(C).$ We have already seen that these
points are exactly the points $v(\Delta_i).$ 
Therefore, by our assumption that the hyperplane 
$H\subset\sp r$
such that $E_1=X\cap H$ does not contain any line of the 
scroll $X$ through a point of $C\cap E_i,$ for any $i=3, 
\ldots ,s,$ it follows that the $P_i$'s are smooth
points for $\Gamma .$ Hence, the points $c_i:=P_i\times a$ 
of $C'\cap F_1$ ($F_1=\sp 2\times a$)
are smooth for both $C'$ and $U.$ Let $T_i$
denote the tangent plane to $U$ at $c_i;$ then $T_i$
contains the tangent line $t$ to $C'$ at $c_i,$
and the line $P_i\times\sp 1.$ We remark that $T_i$
and $F_1$ intersect along the line $T_{\scriptscriptstyle
{\Gamma ,P_i}}\times a.$ Finally, assume that the
line $L$ on $F_1$ does not contain any point $c_i.$
Then, the linear spans $\langle t\cup L\rangle$ and
$\langle (T_{\scriptscriptstyle{\Gamma ,P_i}}\times a)
\cup L\rangle$ are the same. Now, it is easily seen
by a direct computation that, if the lines $R$, $S$
of the scroll $\sp 2\times\sp 1$ are different, then
for any plane $M$ on $\sp 2\times\sp 1$ the linear
spans $\langle R\cup M\rangle$ and $\langle S\cup M
\rangle$ are also different. From all this it follows
that at $P$ the curve $B_2$ has $deg \ 
E_{\scriptscriptstyle {\infty}} - d$
{\sl different} tangent lines.

We can consider $\pi_L:C'\to B_2$ as the normalization
of the curve $B_2,$ since $C'$ is smooth and $\pi_L$
is birational. Moreover, the map $\pi_L$ is unramified
at every point $c_i;$ in fact, in $\sp 3=\langle
t\cup L\rangle$, where $t$ denotes the tangent line to $C'$ 
at $c_i,$ the planes $\langle P\cup t\rangle$ and 
$\langle P\cup L\rangle$ are different. From this it
follows that {\sl the branches of $B_2$ at $P$ are nonsingular}
([5]).

Since $\alpha$ induces a birational map between 
the curves $C$ and $C'$ and since $\pi_L$ does the 
same  between the curves $C'$ and $B_2,$ the geometric 
genus of $B_2$ is equal to the geometric genus of $C,$ 
namely to the sectional genus of $X.$

To complete the proof of (iii) it remains to show
that $B_2$ is smooth outside $P.$ We start by remarking 
that {\sl the cone $\pi_L(U)=\Phi\subset\sp 3$ depends
only on the choice of the plane $F_1$ on 
$\sp 2\times\sp 1,$ and not on the center of projection
$L$ on $F_1.$} In fact, let $L'\subset F_1$ be
another line. Then, for the general line $b\times\sp 1$ on 
the ruled surface $U$ we have that the linear span of
$L$ and $b\times\sp 1$ coincides with the linear span 
of $L'$ and $b\times\sp 1.$ Therefore, on $\Phi$ we
have the family of curves $\{\pi_L(C')\}$ with $L$ varying
in $\check F_1.$ Moreover, it is clear that all these
curves form a linear system $\Omega$ on $\Phi .$

Now, $C'$ intersects $F_1$ in at least two points 
because, otherwise, $C'$ would be birationally
isomorphic to $\sp 1,$ hence rational.
But this contradicts our assumption that
the sectional genus of $X$ is $>0.$ So, we can choose
as $L$ a chord of $C'$ contained in $F_1.$ 

{\sl The curve $\pi_L(C')$ is smooth outside
$P.$} In fact, a singular point would correspond to a
chord $M$ of $C'$ intersecting  $L.$ In this case
the plane $\langle M\cup L\rangle$ has four
points in common with $\sp 2\times\sp 1,$ hence it is
contained in $\sp 2\times\sp 1,$ since $deg(\sp 2\times\sp 1)
=3.$ Then $\langle M\cup L\rangle =F_1$ because 
the only planes on $\sp 2\times\sp 1$ are those of the
ruling, and we are done.

\noindent
Therefore, for the general line $L\subset F_1$ the curve
$\pi_L(C')$ is smooth outside $P.$

\smallskip

By Proposition 2.2, to get the surface of 
$\vert\Sigma\vert$ representing the hyperplane section 
$S$ of $X$ with a hyperplane through $\Lambda ,$
we have just to add the cone $\Phi$ to the plane $g(S).$
Hence, the degree of the surfaces $\Sigma$ is $n+1.$
 
Let $S$ be an arbitrary hyperplane section of $X.$ 
Note that $S\cap E_1\neq\emptyset ,$ hence $P\in\Sigma =
g(S).$ Since the lines $R\subset\sp 3$ through $P$ 
correspond to the lines of the scroll, a general $R$ 
intersects a fixed surface $\Sigma$ outside $P$ in 
only one point and $\Sigma$ is a monoid 
with $P$ as $n$-ple point. 
The assertion concerning the tangent cone of $\Sigma$ at 
$P$ is easily checked and the proof of (iv) is complete.

\smallskip

Consider the restriction of $g$ to $E_i,$ namely
$g:E_i\to B_i;$ if we restrict further on to a hyperplane
section of $E_i$, we get a covering of $B_i$ 
of degree $\delta_i:$ this proves the last assertion.
\endproof

\medskip

\noindent
{\bf Remark 3.2} As observed in the course of the above
proof, the sections $S_a$ of $X$ with the hyperplanes 
$H_a$ of $\sp r$ containing $\Lambda$ correspond to the 
surfaces of $\vert\Sigma\vert$ which break into the cone
$\Phi$ and a plane $W_a$ containing the line $B_1$
(see also Remark 2.3). The surfaces $\Sigma\in
\vert\Sigma\vert$ cut $W_a$ along curves of degree $n+1$
all having $B_1$ as a component. By eliminating this
fixed component we get a linear system $\vert\Sigma '\vert$
of plane curves of degree $n$ on $W_a.$ This linear
system defines the inverse of the birational map $g:
S_a\dashrightarrow  W_a=\sp 2.$ In the case $\vert\Sigma 
'\vert$ is well understood and unique (see, e.g. [1],[2]), 
we get useful informations about $\vert\Sigma\vert.$

We will show how this remark can be used in the last
section, where we will deduce from it the completeness
of the linear system $\vert\Sigma\vert$ constructed above
in the case of the Palatini scroll.

\medskip

\noindent
{\bf Remark 3.3} The restriction of $h:X\dashrightarrow
\sp 2$ to the hyperplane section $E_1$ of $X$ is a
birational map $h:E_1\dashrightarrow\sp 2,$ whose
inverse is given by a linear system $\vert nH_{
\scriptscriptstyle{\sp 2}}-Z\vert,$ where $Z$ denotes
the base locus.

{\sl If the general hyperplane section $S$ of $X$ is
linearly normal (i.e.: if the linear system cut out on
$S$ by the hyperplanes of $\sp r$ is complete), 
then $h(D)\in\vert nH_{\scriptscriptstyle{\sp 2}}-Z
\vert$ for every curve section $D$ of $X.$}

\smallskip

In fact, let $H$ be a hyperplane of $\sp r$ containing 
$D$ and set $S:=H\cap X$ \  and \  $D':=H\cap E_1.$   
On $S$ the curves $D$ and $D'$ are linearly equivalent,
hence $h(D)$ and $h(D')$ are linearly equivalent on $\sp 
2.$ From this it follows that the surfaces $h^{-1}h(D)$ 
and $h^{-1}h(D')$ are linearly equivalent on $X,$ and
therefore $D'=E_1\cap h^{-1}h(D')$ is linearly equivalent 
to $E_1\cap h^{-1}h(D)$ on $E_1.$ Since $D'$ is a hyperplane 
section of $E_1,$ by our assumption $E_1\cap h^{-1}h(D)$ 
is also a hyperplane section of $E_1,$ and we can conclude
$h(D)=h(E_1\cap h^{-1}h(D))\in\vert nH_{
\scriptscriptstyle{\sp 2}}-Z\vert.$


\tit
{\bf 4.-\ Examples: the scrolls of Bordiga and Palatini.}
\tit 


In this section we will apply Theorem 3.1 to the scrolls of
Bordiga and Palatini in $\sp 5.$

\medskip

Let $X\subset\sp 5$ be a Bordiga scroll. It has degree $6$, 
sectional genus $3$  and it a scroll over $\sp 2,$ the scroll 
map $u :X\to \sp 2$ being the adjunction map.
The hyperplane sections of $X$ are Bordiga surfaces; 
in particular, $S\subset\sp 4$ is the image 
of a rational map (birational onto $S$) $ \phi 
:\sp 2\dashrightarrow \sp 4$ defined by the 
linear system $\vert 4H_{\scriptscriptstyle {\sp 2}}
-\sum_{\scriptscriptstyle{1\leqslant i\leqslant 10}} 
y_i\vert$ ([1]). The inverse of $\phi$ is defined 
on the whole $S$ and is the adjunction map.

Therefore, for the Bordiga scroll the map $h$ is 
equal to the scroll map $u.$ In particular, there 
are no exceptional divisors $\Delta$ for the map 
$v:V\to\sp 2.$ Hence, the only exceptional divisors
for $g:X\dashrightarrow\sp 3$ are $E_1,$ $E_2,$ and 
$E_{\scriptscriptstyle {\infty}},$ and the irreducible 
components of $Bs(\vert\Sigma\vert )$ are $P,$ 
$B_1,$ and $B_2.$

The surfaces $\Sigma$ have degree $5.$ Since the scroll 
map $u :X\to\sp 2$ is the adjunction map and 
$deg\thinspace\Gamma =4,$ we can compute the degree of 
$E_{\scriptscriptstyle {\infty}}$ as follows:

\vskip -0.2 cm

$$
deg \ E_{\scriptscriptstyle {\infty}}=E_{\scriptscriptstyle 
{\infty}}\cdot H_X^2=u^{-1}(4H_{\scriptscriptstyle 
{\sp 2}})\cdot H_X^2=4(2H_X+K_X)\cdot H_X^2
$$

\smallskip

\noindent
But $K_X\cdot H_X^2=-8$ ([7]), and we conclude $deg \ 
E_{\scriptscriptstyle {\infty}}=16.$  Therefore, the degree 
of $B_2$ is $14,$ and the point $P$ has multiplicity $10$ 
for $B_2.$

We recall that the exceptional divisor $E_1$ for $g$ is a 
hyperplane section of $X.$ 
On $E_1$ there are exactly ten lines of the scroll, and the 
adjunction map $h$ contract these lines to the points
$P_1,\ldots ,P_{10}\in\sp 2.$ By [4], Prop. 
4.1, we can assume that these points satisfy any kind
of \lq\lq\thinspace general position" assumption.
Moreover, observe that, by construction of the map $g,$
the set $\{ P_1,\ldots ,P_{10}\}\subset\sp 2$ is isomorphic
to the plane section of the tangent cone to $B_2$ at $P.$

The first conclusion we can draw from these remarks 
is that {\sl a
surface of degree $5$ in $\sp 3$ containing $B_2$ is
necessarily a monoid of vertex $P$}
(consider the Taylor expansion at $P$ of the equation
of the surface). 

In view of computing the arithmetic genus of $B_2$
starting from its geometric genus $p_{\scriptscriptstyle g}
(B_2)=3,$ it is advisable
to have more precise informations on the structure
of the singular point $P.$ We let $A$ denote the local
ring of $B_2$ at $P,$ $\cal M$ its maximal ideal, and, 
finally, we let $G$ denote the associated graded ring.
{\sl $G$ is reduced.} In fact, we can apply again the
above remarks and assume that the points $P_1,\ldots 
,P_{10}\in\sp 2$ are in \lq\lq\thinspace generic position"
(definition in [6]). Then, the conclusion follows 
from $(iii)$ of Thm.2.4. and [6], Thm. 4.2.

Now, let $\overline A$ denote the integral closure of
$A$ into its field of fractions. From [6], Thm. 4.4
we conclude that ${\cal M}^3$ is the conductor of $A$
into $\overline A.$ Therefore

$$dim \biggl( {A\over{{\cal M}^3}}\biggr) =10
\ \ \ \ \ \ \ \ \ \ \ 
\hbox{and}\ \ \ \ \ \ \ \ \ \ \ 
dim \biggl({{\overline A}\over{{\cal M}^3}}\biggr)
=10\cdot 3=30$$

\smallskip\noindent where $10$ is the multiplicity
of the point $P$ on $B_2.$ Hence $\delta =dim({\overline A}
/A)=20,$ and we conclude

$$
\ga (B_2)=p_{\scriptscriptstyle g}(B_2)+\delta =23.
$$

\medskip
\noindent
{\bf Remark 4.1}\ In [4] we got a different linear system
$\vert\Sigma'\vert$ of quintic surfaces in $\sp 3$
defining a birational map $f':\sp 3\dashrightarrow X$.
The inverse of $f'$ is the projection centered at a suitable
line $L'\subset X$. 
If we compose $f':\sp 3\dashrightarrow X$ with $g:X\dashrightarrow\sp 3$
we clearly get an element $T$ of the Cremona group  of $\sp 3$,
which depends on the choices of $L'$ for $f'$ and
of $\Lambda$ and $L$ for $g$.

For general choices, $T$ is of type $(7,5)$, i.e.
the surfaces defining $T$ have degree $7$
and those defining $T^{-1}$ have degree $5$.
But, for special choices, this degrees can decrease.
In particular, if $L'\subset\Lambda$, we get 
a cubo--cubic transformation.

Let us analyze briefly the situation in the case
$\Lambda\cap L'=\emptyset$.

$T^{-1}$ is defined by the subsystem of $\vert\Sigma\vert$
of surfaces containing $g(L')$. It is easy to
prove that this curve is a skew cubic. Let $M$ be a
plane in $\sp 3$: the surface $f(M)$ has degree $5^2-\deg(B_1
\cup B_2)=10$. Since $\sharp(f(M)\cap L')=\sharp(M\cap g(L')=3$, 
$T(M)$ is a surface of degree $7$.

If $L'\subset\Lambda$, $L'$ is a component of $C=X\cap\Lambda$,
hence it is part of the fundamental locus of $g$.
The adjunction map $u:X\to\sp 2$ takes both $L'$ and its
residual in $C$ to conics. Hence the surface $\Phi$ (Prop. 2.2)
splits in two quadrics cones. One of them is contained in all
surfaces $\Sigma$ corresponding to hyperplane sections
of $X$ containing $L'$. This shows that $T^{-1}$ is
defined by cubics. Conversely, it is rather easy to see that
the image of a line via $T$ is a skew cubic ([4], Prop. 3.11), so also
$T$ is defined by cubic surfaces.
\medskip
Let $X\subset\sp 5$ be a Palatini scroll; 
it has degree $7$ and sectional genus $4.$  
A Palatini scroll is a scroll over a smooth cubic surface 
$V\subset\sp 3,$ and the scroll map $u :X\to V$ is the 
adjunction map. If $S$ denotes a general hyperplane section 
of $X,$ then the restriction of $u$ to $S$
is the adjunction map for $S.$ 
The natural choice for the map $v:V\to\sp 2$ in our case
is the blow-up of $\sp 2$ at the points 
$x_1,\ldots ,x_6\in\sp 2.$ The composition $S\to\sp 2$ 
of these maps is the blow-up of 
$\sp 2$ at  eleven points  $x_1,\ldots ,x_6,y_1,\ldots 
,y_5;$ more precisely, the linear system of curves in 
$\sp 2$ defining the birational inverse of $S\to\sp 2$ 
is $\vert 6H_{\scriptscriptstyle {\sp 2}}-\sum_{
\scriptscriptstyle{1\leqslant i\leqslant 6}}2x_i
-\sum_{\scriptscriptstyle{1\leqslant j\leqslant 5}} 
y_j\vert$ (see [1]). Therefore, for the Palatini scroll 
we have $n=6$ and the surfaces $\Sigma$ will have degree $7.$ 

A computation similar to that performed above for the 
Bordiga scroll shows that $deg \ 
E_{\scriptscriptstyle {\infty}}=12.$ Hence, the degree of $B_2$ 
is $11$ and the multiplicity of $P$ for $B_2$ is $5.$

The exceptional divisors for the map $v:V\to\sp 2$ are six
lines on $V.$ It is shown in [4] that, for any line $W$ 
on $V,$ the surface $u^{-1}(W)$ is a quadric on $X.$
Therefore, as components of the base locus of $\vert\Sigma\vert$ 
we have also the first infinitesimal neighborhoods of
six lines $B_3,\ldots ,B_9,$ each containing $P.$ 
Let us denote by $D$ the first infinitesimal 
neighborhoods of $B_3\cup\ldots\cup B_9.$ 
Note that for any surface in $\sp 3$ containing $D$ the
point $P$ has multiplicity at least $6.$ In particular,
any surface of degree $7$ containing $D$ is already a monoid.

The curve  $B_2$ intersects the line $B_1$ at $6$ points
and  the lines $B_3,\ldots ,B_9$ are chords for $B_2$.

The geometric genus of $B_2$ is $4;$ it is possible to
perform a computation similar to that for Bordiga to
show that $\ga (B_2)=10.$ From this one can compute
$\ga (Bs(\vert\Sigma\vert ))=97.$

\medskip

We want to prove now that 

\smallskip 

\noindent
{\bf Proposition 4.2}\ {\sl The linear system
$\vert\Sigma\vert$ constructed above to represent the
Palatini scroll is complete.}

\proof The cone
$\Psi$ projecting $B_2$ from $P$ has degree $11-5=6,$
where $5$ is the multiplicity of $P$ on $B_2.$
Moreover, $\Psi$ is clearly singular along any line $B_3,
\ldots ,B_9,$ hence $\Psi$ contains $D.$ As remarked above,
for any surface containing $D$ the
point $P$ has multiplicity at least $6.$ Then, a
surface of degree $6$ containing $D$ is necessarily a
cone with vertex $P,$ and we conclude that {\sl $\Psi$ is
the only surface of degree $6$ containing $D\cup B_2.$}

Since $E_1$ contains $C=X\cap\Lambda ,$ and is a hyperplane
section of $X$ it corresponds to a reducible surface 
$\Sigma$ having as a component a plane $W$ containing the 
line $B_1$ (see Remark 3.2).
Moreover, the hyperplane sections of $E_1$ correspond to
the curves of degree $6$ of a suitable linear system 
$\vert\Sigma '\vert$ on $W.$ We denote by $B$ and $Z$
the base loci of $\vert\Sigma\vert$ and $\vert\Sigma '
\vert$ respectively. 

Then we can define a linear map

$$
\varphi : H^0(\sp 3,\hbox{\fascio I}_B(7))\to 
H^0(W,\hbox{\fascio I}_Z(6))
$$

\smallskip\noindent by setting \  $\varphi (F):=\overline F
\cdot h^{-1},$ where $\overline F$ denotes the restriction
of $F$ to $W,$ and $h$ denotes the equation on $W$ of the
curve $W\cap\Psi .$ The kernel of $\varphi$ is the one
dimensional vector space generated by the equation of $W
\cup\Psi .$ 

We know that the linear system $\vert 6H_W-
Z\vert$ is complete and of dimension $4,$ namely
$h^0(W,{\fascio I}_Z(6))=5$ ([2]). Therefore
$h^0(\sp 3,{\fascio I}_B(7))=6$ and the proof is complete.
\endproof

\medskip

We conclude with the following remark.

\smallskip

\noindent
{\bf Remark 4.3}\ {\sl The scrolls of Bordiga and
Palatini are both suitable projections into $\sp 5$ of the
Veronesean embedding of $\sp 3$ by the linear system
of surfaces of degree $7.$}

Let $X$ be a Bordiga scroll.  By Theorem 3.1 it is 
sufficient to construct a rational map $h:X
\dasharrow\sp 2$ such that for the general curve section 
$D$ of $X$ the plane curve $h(D)$ has degree $6.$

\smallskip

Since the general hyperplane 
section $S\subset\sp 4$ of $X$ is linearly normal,
by Remark 3.3  we have $h(D)
\in\vert 4H_{\scriptscriptstyle{\sp 2}}-\sum_{
\scriptscriptstyle{1\leqslant i\leqslant 10}} x_i\
\vert$ for every curve section $D$ of $X.$

\smallskip

Let $w:\sp 2\dasharrow\sp 2$ denote the standard quadratic 
transformation centered at $x_1,x_2,x_3,$ and let $y_1,y_2,
y_3$ denote the points in which the three lines $x_ix_j$ are 
contracted by $w.$ If we compose
$w$ with $u$ we get a new rational map $X\dasharrow\sp 2$
that we will denote again by $h,$ such that
for every curve section $D$ of $X$ we have $h(D)
\in\vert 5H_{\scriptscriptstyle{\sp 2}}-2\sum_{
\scriptscriptstyle{1\leqslant i\leqslant 3}} y_i-\sum_{
\scriptscriptstyle{4\leqslant i\leqslant 10}} x_i\
\vert .$

\smallskip

We repeat this procedure by composing the rational map 
$X\dasharrow\sp 2$ obtained above with the standard 
quadratic transformation $t:\sp 2\dasharrow\sp 2$ 
centered at $y_3,x_4,x_5.$ For the new rational map 
$h:X\dasharrow\sp 2$ we have $h(D)\in\vert 
6H_{\scriptscriptstyle{\sp 2}}-3z-2\sum_{
\scriptscriptstyle{1\leqslant i\leqslant 4}} y_i-\sum_{
\scriptscriptstyle{6\leqslant i\leqslant 10}} x_i\
\vert$ (with a suitable labelling of the points),
for every curve section $D$ of $X.$


\tit
\centerline{\bf References}
\tit


\itemitem {1.} {\caps J. Alexander,} 
Surfaces rationnelles non-sp\'eciales dans $\sp 4$\hfill\break
{\it Math. Z.\/} {\bf 200} (1988), 87--110.

\vskip 0.1 cm 

\itemitem {2.} {\caps F. Catanese and M. Franciosi,}
Divisors of small genus on surfaces and projective embeddings 
\hfill\break
{\it Israel Mathematical Conference Proceedings\/}  {\bf 9} (1996),
109--140.

\vskip 0.1 cm 

\itemitem {3.}  {\caps F. Jongmans,}
Les vari\'et\'es alg\'ebriques \` a trois dimensions dont les 
courbes - sections ont le genre trois\hfill\break 
{\it Acad. Roy. Belgique, Bull. Cl. Sci.\/} (5) {\bf 30} (1943), 
766--782 and 823--835.
 
\vskip 0.1 cm

\itemitem {4.} {\caps E.  Mezzetti and D. Portelli,} 
On smooth rational threefolds of $\sp 5$ with rational 
non-special hyperplane section\hfill\break 
preprint (1997), to appear in {\it Mathematische Nachrichten\/}.

\vskip 0.1 cm

\itemitem {5.} {\caps  F. Orecchia,}
Ordinary singularities of algebraic curves\hfill\break
{\it Canad. Math. Bull.\/}
{\bf 24} (4) (1981),  423--431.
 
\vskip 0.1 cm 

\itemitem {6.} {\caps  F. Orecchia,}
Points in generic position and conductors of
curves with ordinary singularities\hfill\break
{\it J. London Math. Soc.\/}
{\bf 24} (1981),  85--96.
 
\vskip 0.1 cm 

\itemitem {7.} {\caps  G. Ottaviani,}
On $3$-folds in $\sp 5$ which are scrolls\hfill\break
{\it Ann. Scuola Norm. Sup. Pisa, Cl.
Sci.\/} Ser. (4) {\bf 19} (1992),  451--471.

\bye